\documentclass{amsart}
\usepackage{amsmath,amsthm,amsfonts, amssymb,amscd}
\usepackage[all]{xy}
\usepackage{mathrsfs}
\newtheorem{theorem}{Theorem}[section]
\newtheorem{lemma}[theorem]{Lemma}
\newtheorem{definition}[theorem]{Definition}
\newtheorem{proposition}[theorem]{Proposition}

\begin{document}
\title{Associativity of $\sigma$-sets for non-antielement $\sigma$-set groups}
\author{Alfonso Bustamante Valenzuela}

\address{Department of Informatics, Technological University of Chile (INACAP), Apoquindo 7282, Santiago, Chile}
\email{alfonso.bustamante02@inacapmail.cl}

\date{}%
\maketitle

%

\begin{abstract}
We study and extend the conditions for asociativity on fusion over antielement free $\sigma$-sets to introduce a group to solve $\sigma$-set equations.  $\sigma$-sets as a theory of sets and antisets is sumarized and used as a framework to define the main elements of this work. A theorem on Local Asociativity, conmutative groups and solution  of one variable fusion equation is presented.  
\end{abstract}

\section{Introduction}
The inverse element play a very important role, together with associativity, on solving equations form an algebraic structure. But what happend when the inverse element and associativity together arise unsolvable inconsitencies? Let us pose the following example. Let $A$ and $B$ two sets differet to empty, with the propietry that $A \cup B=\emptyset.$ Associativity of a binary operation $\ast$ is defined as $(A \ast B) \ast C=A \ast (B \ast C),$ which in this context would be $(A \cup B) \cup C=A \cup (B \cup C).$ Let us say, finally, that $B=C.$ So we have,
 $$(A \cup B) \cup C=\emptyset \cup C,$$
 $$A \cup (B \cup C)=A \cup B.$$

As $\emptyset \cup C=C,$ where, if associativity holds under our previous assumptions, we would conclude that $B=\emptyset,$ which contradict our first assumption. \\

That is the case of set theory, where this problem arised after the aparition of the concept of antisets. Antisets are a relativley new object studied on set theory, varing from elimination propietry under union \cite{Car09}, \cite{Ar09}, to gaps on a set which eliminate any element that join with \cite{Bal13}; in this article, however, we will strictly use the first type of Antisets, particulary the definition developed by Arraus \cite{Ar09}. \\

What we want to present in this work, is the algebraic condition which aloud local associativity of fussion on $\sigma$-sets, this is, extend the Arraus condition of associativity on antielement free $\sigma$-sets and open the possiblity to associate sets and antisets. Our goal is to to make partition equations of the form $$A \cup X=B,$$ as where presented on \cite{Bus11},  solvabe, for that matter, we will define a special type of group based on the locally associative $\sigma$-set. \\

Before we present our results, we want to introduce some notation of Arraus work, which will help us to analize and define some of the concepts needed to proof our main theorem on associativity on $\sigma$-sets.\\

\section{Preeliminars on $\sigma$-sets}

In \cite{Ar09} it was defined a type of set denoted as $\sigma$-set, which was the gathering of both sets and antisets, where for a set $A$ and a set $A^{\ast}$ on a $\sigma$-set, the following identity holds,
$$A \cup A^{\ast}=\emptyset,$$
which is exactly replied on their elements, this is, for $X=\{x\}$ and $X^{\ast}=\{x^{\ast}\},$ $\{x,x^{\ast}\}=\emptyset.$ 

For X and Y $\sigma$-sets, he defined two base operations,

\begin{enumerate}
\item[{\rm (1)}] $X \hat{\cap} Y:=\{x \in X:x^{\ast} \in Y\}$
\item[{\rm (2)}] $X \divideontimes Y:=X-X \hat{\cap} Y.$
\end{enumerate}

Below, we paraphrase an example appling the above definitions.

Let $X=\{1,2\}$ and $Y=\{1^{\ast},2^{\ast} \},$ where,
\begin{enumerate}
\item[{\rm (1)}] $X \hat{\cap} Y=X$
\item[{\rm (2)}] $X \divideontimes Y=\emptyset.$
\end{enumerate}
From this concepts, the author also defined what he denoted as {\bf antielement free} $\sigma$-set, to a set $X,$where for $A,B \in X, A \hat{\cap} B=\emptyset.$ 

This $\sigma$-set was later denoted as,
$$AF=\{x: x \text{is antielement free} \} .$$

A variation of union, called fusion,  presented by Arraus at the begining of the article, is redefined as,
$$X \cup Y=\{x:(x \in X \divideontimes Y) \vee (x \in Y \divideontimes X)\}$$
Finaly, a theorem is introduced in order to show why associativity is not possible in fusion over $\sigma$-sets, this is, on not antielement free $\sigma$-sets.

\begin{theorem} [Th.3.52, p.26]
Let $X$ be a $\sigma$-set. Then
\begin{enumerate}
\item[{\rm (a):}] $X \cup \emptyset=\emptyset \cup X=X$
\item[{\rm (b):}] $X \cup X=X.$
\end{enumerate}	
\end{theorem}

It is shown that fusion is conmutative, this is $X \cup Y=Y \cup X,$ but when associativity is to be shown, the author point out the following case. Let us have the $\sigma$-sets $X=\{1,2\},Y=\{1^{\ast},2^{\ast} \}$ and $Z=\{1\},$ where, 
$$(X \cup Y) \cup Z=\emptyset \cup Z=Z$$
and
$$X \cup (Y \cup Z)=X \cup Y=\emptyset,$$
where, as in the ecample presented in this article, ocurrs the same contradiction where $Z=\emptyset.$

To solve this problem, the author proposed the following theorem,

\begin{theorem}
If $X \in AF,$ then the fusion in $X$ is associative, that is
$$(X \in AF)[(A \cup B) \cup C=A \cup (B \cup C)].$$
\end{theorem}

Under this scheme, the $\sigma$-sets, to have the propietry of being associative, there must be antielement free. But what happend when equations like $A \cup X=B$ has to be solved?  Is there a restriction from which associativity on $\sigma$-sets holds? Such restriction is our main result, a theorem on the condition of associativity on $\sigma$-sets, which we will present in the following section together with some definitions an lemmas to prove it.

\section{Definitions and Main result}
\subsection{Associative chain of fusion and Localy associative $\sigma$-Sets}
Before introduce our main result, let us present the following previous definitions and lemmas.

We are going to take the notion of {\bf chain of fusion} which was defined by Arraus as,
$$\bigcup \overrightarrow{\mathcal{F}}=(A \cup B) \cup C,$$
which for the prupose of the following definitions, we will denote as $\overrightarrow{ABC},$ where $ABC$ represent a fixed order on the expression.

\begin{definition}
	Let $A,B$ and $C$ $\sigma$-sets on $S.$An {\bf Associativity Evaluation Chain} on $S$, denoted as $E_S,$ is an equation on $\sigma$-sets, where,
$$E_{S}=\overrightarrow{ABC} \cup \overrightarrow{C^{\ast}B^{\ast}A^{\ast}}.$$	
\end{definition}

\begin{definition}
An {\bf Associative chain of fusion over S}  is an asociative chain where, for every $A,B$ and $C$ $\sigma$-sets on $S$,
$$(A \cup B) \cup C=A \cup (B \cup C).$$
\end{definition}

\begin{lemma}
Let $A,B$	and $C$ be $\sigma$-sets of $S$. if $(A \cup B) \cup C=A \cup (B \cup C),$ then $\overrightarrow{ABC}=\overrightarrow{CBA}.$
\end{lemma}

\begin{proof}
Is enough to prove that, by conmutativity, $A \cup (B \cup C)=(C \cup B) \cup A$ because $X=A$ and $Y=B \cup C;$ the same with $B \cup C=C \cup B,$ where we have the fusion chain $\overrightarrow{CBA}.$ As we previously assumed that $(A \cup B) \cup C=A \cup (B \cup C),$ and we know by definition that $\overrightarrow{ABC}=(A \cup B) \cup C,$ proving the requested identity.
\end{proof}

\begin{lemma}
If $\overrightarrow{ABC}$ is a $\sigma$-set, then $\overrightarrow{A^{\ast}B^{\ast}C^{\ast}}$ is it  antiset.
\end{lemma}

\begin{proof}
By definition of antiset and chain of fusion over $\overrightarrow{A^{\ast}B^{\ast}C^{\ast}}$, is enoug to show that the fusion of $\overrightarrow{ABC}$ and $\overrightarrow{A^{\ast}B^{\ast}C^{\ast}}$ is actually the empty set.
\end{proof}

\begin{proposition}
Let $A,B$ and $C$ be $\sigma$-sets on $S.$	if $E_{S}=\emptyset,$ then fusion over $S$ is associative.
\end{proposition}
ß
\begin{proof}
	Let us show by contapositive. Assume that $(A \cup B) \cup C \neq A \cup (B \cup C).$ We know that the second expression by the law of conmutativity on fusion, we have fusion chain $ABC.$ If we fusion  $(A \cup B) \cup C,$ this is, the fusion chain $\overrightarrow{ABC}$ with the antielement $\sigma$-set $\overrightarrow{C^{\ast}B^{\ast}A^{\ast}},$ we will obtain,
	$$\overrightarrow{ABC} \cup \overrightarrow{C^{\ast}B^{\ast}A^{\ast}} \neq \emptyset,$$
	because of $\overrightarrow{ABC} \neq \overrightarrow{CBA}$ by our  hypothesis, hence the expression obtained is equivalent, by definition, to the equation $E_{S},$ so the proof is complete.
\end{proof}

Let us ilustrate the above proposition with the following example.

We have the $\sigma$-sets $A=\{a,b\}, B=\{a^{\ast},b^{\ast} \}$ and $C=\{c,d\}.$ Taking the definition of fusion, we have that $$\overrightarrow{ABC}=C$$ and $$\overrightarrow{C^{\ast}B^{\ast}A^{\ast}}=C^{\ast},$$ where finally, the Associative evaluation chain $E_{S}=\emptyset,$ so by the proposition, $ABC$ is an associative chain of fusion.

Now let us take a contrary example to see what happend when the Associative evaluation chain is not $\emptyset.$ We have the $\sigma$-sets $A=\{1,2\}, B=\{1^{\ast},2^{\ast} \}$ and $Z=\{1^{\ast} \}.$ Taking the definition of fusion, we have that $$\overrightarrow{ABC}=Z,$$ and $$\overrightarrow{C^{\ast}B^{\ast}A^{\ast}}=\emptyset,$$ where $$E_{S}=Z \cup \emptyset=Z,$$ which is different to $\emptyset,$ so we can conclude that $ABC$ is not associative. \\

This results as we can see, is a very useful tool to evaluate associativity over an ordered configuration of $\sigma$-sets; however, we are still not able to evaluate on a wide combinatorial range, this is, over the generalized union $\bigcup \{A,B,C\}.$ For thar reason, we develop the following definitions. \\

\begin{definition}
Let $A,B$	and $C$ $\sigma$-sets. $A,B,C$ are said to be {\bf Locally Associative} over fusion if, for any combination of $A,B$ and $C,$ the fusion chain is associative.
\end{definition}

There exist some $\sigma$-sets $S$ such that are associative chain of fusion but are not locally associative. Let us examine the following $\sigma$-set where $A=\{1,2\},B=\{1^{\ast},2^{\ast} \}$ and $C=\{1,2\}.$ As we can see, by definition of associative evaluation chain, we have that $E_{S}=\emptyset,$ so is associative; but if we re-arrange the set on the fusion set as $BAC,$ we will have that the new evaluation chain would be diferent to empty set, impliyng that the associativity is not defined for every configuration.

\begin{definition}
	For the equation system $E=\{E_{X},E_{Y},E_{Z}\},$ where the $\sigma$-sets $X,Y,Z$ are in $E,$ the equations 
$$E_{X}=\overrightarrow{XYZ} \cup \overrightarrow{Z^{\ast}Y^{\ast}X^{\ast}},$$
$$E_{Y}=\overrightarrow{YZX} \cup \overrightarrow{X^{\ast}Z^{\ast}Y^{\ast}},$$
$$E_{Z}=\overrightarrow{ZXY} \cup \overrightarrow{Y^{\ast}X^{\ast}Z^{\ast}}.$$
\end{definition}



	


\begin{theorem}
	For $X,Y$ and $Z$ $\sigma$-sets, if $E= \{\emptyset \},$ then fusion is locally associative on $\sigma$-sets.
\end{theorem}

\begin{proof}
As $E=\{\emptyset\}$ we know that the equations $E_{X},E_{Y}$ and $E_{Z}$ are equal to the empty set, wich implies that $XYZ,YZX$ and $ZXY$ are conmutative fusion chains. Because of lemma 3.3, same $ZYX,XZY$ and $YXZ$ are. As every configuration is an associative fusion chain, by definition, is locally associative.
\end{proof}

\subsection{Main Theorem on  $\sigma$-sets Conmutative groups}

To start this section, let us remind what is a group as an algebraic structure. First, we remember that an algebraic structure \cite{Bo70} is a pair $(A,\ast)$ where $A$ is a set of certain mathematical objects, and $\ast$ is an action from $A \times A$ to $A,$ this is, the result of the operation of two elements of $A$ also belong to $A.$ \\

Algebraic structures posess also composition laws, this is, laws of what happend when operations do certain combination of special objects of the set or coclatenations like associativity. 

A group $(G,\ast)$ is a set $G$ of elements a,b,c such that,
\begin{enumerate}
\item[{\rm (G1)}] $a \ast e=a$
\item[{\rm (G2)}] $a \ast a^{-1}=e$
\item[{\rm (G3)}] $(a \ast b) \ast c=(a \ast b) \ast c,$
\end{enumerate}
if in the group $a \ast b=b \ast a,$ it said to be conmutative or abelian. The group we are going to study, is an abelian one, we will proof that the theorem 3.8. alows the minumal conditions to define a cancelative group.


\begin{theorem}
	Let $A,B,C$ be $\sigma$-sets on $\mathbb{G}.$ if $\mathbb{G}$ is locally associative, then $(\mathbb{G},\cup)$ is an associative group.
\end{theorem}

\begin{proof}
As $\mathbb{G}$	is locally associative, means that $E=\{\emptyset \}$ so, there is neutral element; in the other hand, as every $E_{X},E_{Y},E_{Z}$ in $E$ has in their expression an antiset, where by definition of evaluation chain and lemma 3.2, the rest of the needed condition to be a group are acomplished.
\end{proof}

\begin{theorem}
	If $(S,\cup)$ is a Locally associative group, then $X \cup A=B$ is solvable and the solution is $X=B \cup A^{\ast}.$
\end{theorem}

\begin{proof}
the $\sigma$-set $A,B,C$ is locally asociative, so, if we left add to both sides of the equation the term $A^{\ast},$ the right expresion $(X \cup A) \cup A^{\ast},$ by local associativity and cancellation with antisets, you will have $X,$ the left term remains without change so we can conclude that has a solution. 
\end{proof}






To have a better understandig of this result, let us ilustrate it with the following example. 

Let us have a $\sigma$-set $S$ where $A=\{\alpha,\beta\},$ and $B=\{a^{\ast},b^{\ast},c^{\ast},\alpha,\beta\},$ and X is to be found in the form
$A \cup X=B.$
First, we separate the expression
$$\{a^{\ast},b^{\ast},c^{\ast},\alpha,\beta\}=\{a^{\ast},b^{\ast},c^{\ast}\} \cup \{\alpha,\beta\}$$
then, by definition of locally associative $\sigma$-set, we know that $X$ is part of an associative chain of fusion over any permutation. This aloud us to just aplly the theorem on result for this kind equations, so we have,
$$X=(\{a^{\ast},b^{\ast},c^{\ast}\} \cup \{\alpha,\beta\})\cup \{\alpha^{\ast},\beta^{\ast} \},$$
because this is an associative chain, is enough to apply associativity and cancelation over $\sigma$-sets to finaly obtain that $X=\{a^{\ast},b^{\ast},c^{\ast}\}$ which is the solution to obtain the sigma set $B.$
\section{Further applications}

A possible application of this research is on link algebra, where by the theorem on locally associative $\sigma$-sets over fusion, it is posible to make graph equations over union, as equations of the form $A \cup X=B,$ as where presented on \cite{Bus11},  solvabe. \\

There is an other case where this kind of special $\sigma$-set space could be appliyed: graph join equations, where a graph $G_{1}$ is joined to a graph $G_{2}$ if al the vertices of the two graphs are joined by edges. This type of result, could constitute a case of a join graph equation of one variable, where the inverse graph $G_{2}^{\ast}$ could provide the graph to obtain the resultant graph, this would be specially useful on Zykov's Linear Complexes \cite{Zy49}, where graph join was denoted as $G_{1}G_{2},$ in order to solve problems on graph concentration. 


\end{document}